\documentclass[10pt]{article}
\usepackage{amsgen, amsmath, amsfonts, amsthm, amssymb, enumerate,amscd}

\newtheorem{defn}{Definition}[section]
\newtheorem{prop}[defn]{Proposition}
\newtheorem{lem}[defn]{Lemma}
\newtheorem{thm}[defn]{Theorem}

\newtheorem{rem}[defn]{Remark}

\newcommand {\ZZ}{{\mathbb Z}}

\newcommand {\XX}{{\mathcal X}}
\newcommand {\K}{{\mathcal K}}
\newcommand {\E}{{\mathcal E}}

\newcommand {\G}{{\mathbb  G}}
\newcommand {\Z}{{\mathcal Z}}
\newcommand {\CC}{{\mathcal C}}
\newcommand {\Q}{{\mathbb Q}}
\newcommand {\R}{{\mathbb R}}
\newcommand {\OO}{{\mathcal O}}

\newcommand {\M}{{\mathcal M}}

\newcommand {\CP}{{\mathbb P}}

\newcommand {\D}{{\mathcal D}}
\newcommand {\QL}{{{\mathbb Q}_{\ell}}}

\newcommand {\J}{{\mathcal J}}

\newcommand {\mE}{{\bar{\mathcal E}}}
\newcommand {\mF}{{\mathcal F}}

\def\div{\operatorname{div}}

\def\dim{\operatorname{dim}}

\def\ord{\operatorname{ord}}

\title{A non-Archimedean analogue of the Hodge-$\D$-conjecture for products of elliptic curves}
\author{Ramesh Sreekantan (Tata Insitute of Fundamental Research)}

\begin{document}

\maketitle
\begin{abstract} In this paper we show that the map
$$\partial:CH^2(E_1 \times E_2,1)\otimes \Q \longrightarrow PCH^1(\XX_v)$$
is surjective, where $E_1$ and $E_2$ are two non-isogenous
semistable elliptic curves over a local field, $CH^2(E_1 \times
E_2,1)$ is one of Bloch's higher Chow groups and  $PCH^1(\XX_v)$ is
a certain subquotient of a Chow group of the special fibre $\XX_{v}$
of a semi-stable model $\XX$ of $E_1 \times E_2$. On one hand, this
can be viewed as a non-Archimedean analogue of the
Hodge-$\D$-conjecture of Beilinson - which is known to be true in
this case by the work of Chen and Lewis \cite{lech}, and on the
other, an analogue of the works of Spei{\ss} \cite{spie}, Mildenhall
\cite{mild} and Flach \cite{flac} in the case when the elliptic
curves have split multiplicative reduction.

\end{abstract}

\section{Introduction}

The aim of this paper is to prove a special case of the following
conjecture: Let $K$ be a local field with residue characteristic $v$
and ring of integers $\OO$. Let $X$ be a variety over $K$ and let
$\XX$ be a semi-stable model over $\OO$. Then the map
$$CH^a(X,b)\otimes \Q \stackrel{\partial}{\rightarrow} PCH^{a-1}(\XX_{v},b-1)$$
is surjective. Here, assuming the Parshin-Soul\'{e} conjecture
\cite{soul}, if $b>1$, $PCH^{a-1}(\XX_v,b-1)$ is the higher Chow
group $CH^{a-1}(\XX_v,b-1)\otimes \Q$. In particular, it is $0$ if
$\XX_v$ is non-singular. If $b=1$ it is a certain subquotient of the
Chow group $CH^{a-1}(\XX_v)\otimes \Q$.

A conjecture of Bloch's and the work of Consani \cite{cons} suggests
that the dimension of the group $PCH^{a-1}(\XX_v,b-1)$ is the order
of the pole of $L_v(H^{2a-b-1}(X),s)$ at $s=(a-b)$. In general this
can be non-zero so the surjectivity is non-trivial.

We prove this in the case when $X=E_1 \times E_2$ where $E_1$ and
$E_2$ are non-isogenous elliptic curves over $K$, $a=2$ and $b=1$.
There are several subcases to be considered, depending on whether
the reduction is good or bad.

\begin{itemize}

\item When $v$ is a prime of good reduction, the expected dimension of
$PCH^{1}(\XX_v,0)$ is $6$,$4$ or $2$, depending on whether the
special fibres $\E_{1,v}$ and $\E_{2,v}$ are isogenous or not. The
surjectivity of the map in this case was shown by Spie\ss
\cite{spie}.

\item When $v$ is a prime of good reduction for one of the elliptic
curves and bad semi-stable reduction for the other the expected
dimension is $2$. In this case it is easy to see what the elements
of $CH^2(E_1 \times E_2,1)$ are.

\item When $v$ is a prime of bad semi-stable reduction for both $E_1$
and $E_2$. In this case the expected dimension is $3$. It is easy to
find elements of $CH^2(X,1)$ which whose image under the map
$\partial$ spans two dimensions. To find an element whose boundary
spans the third dimension seems to require more work, and is the
purpose of this paper.

\end{itemize}

Beilinson's Hodge-$\D$-conjecture \cite{jann},  specialized to our
case, states that the map
$$r_{\D}\otimes \R:CH^2(X,1)\otimes \R \longrightarrow H^3_{\D}(X,\R(2))$$
is surjective. This is now a theorem of Chen and Lewis \cite{lech}.
As explained below, the group $PCH^1(Y)$ shares many properties with
the Deligne cohomology group, so our statement can be viewed as a
non-Archimedean analogue of this. In general the Hodge-$\D$
conjecture is false \cite{mull}.

An $S$-integral version of the Beilinson conjectures, or a special
case of the Tamagawa number conjecture of Bloch and Kato, would
assert that, for a variety over a number field, there are elements
in the higher Chow group over the number field itself which bound
the elements of the Chow groups of the special fibres.

The only case for which there is some evidence is the work of Bloch
and Grayson \cite{blgr} on $CH^2(E,2)$ of elliptic curves, but even
here, as far as I am  aware, there is not a single case of an
elliptic curve over $\Q$ where it is known that there are as many
elements of the group as would be predicted by the full $S$-integral
Beilinson conjecture.

On the other hand, Beilinson  \cite{beil} proved his conjecture for
the product of two non-isogenous elliptic curves over $\Q$, and this
can be viewed as a statement for the Archimedean prime. One expects
\cite{mani}  that the Archimedean prime behaves like a prime of
semi-stable reduction. Further, the work of Mildenhall \cite{mild}
provides evidence for this conjecture when $v$ is a prime of good
reduction and $E_1$ and $E_2$ are isogenous elliptic curves over
$\Q$.

One might hope that one can extend these results to the case when
$E_1$ and $E_2$ are not isogenous and $v$ is a prime of bad
semi-stable ( that is, split multiplicative ) reduction for both of
them. Towards that end, we considered this local situation first.

The outline of the paper is as follows. In the first section we
define the group $PCH$ that appears as the target of the boundary
map. We then specialize to the product of two semi-stable elliptic
curves and describe the fibre of the semi-stable model of the
product of the two curves and the group $PCH$ in this case. Then we
use  the work of Frey and Kani on the existence of curves of genus
$2$ on products of elliptic curves along with Spei{\ss}'s work to
construct some elements in the higher Chow group. Finally we compute
their boundary and show that they suffice to prove surjectivity.

The method of proof is almost identical to that of Spie{\ss}, the
only difference being that we have to modify his arguments
appropriately to work in the case of bad reduction. He obtains some
consequences for codimension $2$ cycles on $\XX$ which follow from
the surjectivity of $\partial$, so they apply in our case as well.

I would like to thank S. Kondo, C. Consani, B. Noohi, C-G Lehr and
C.S. Rajan for some useful conversations, and K. Kimura for pointing
out some mistakes in an earlier version of this paper. I would also
like to thank the Max-Planck-Institut f\"{u}r Mathematik for
providing a wonderful atmosphere in which to work in.

\section{Preliminaries}

Let $X$ be a smooth proper variety over a local field $K$ and $\OO$
the ring of integers of  $K$ with closed point $v$ and generic point
$\eta$.

By  a model $\XX$ of $X$ we mean a flat proper scheme $\XX
\rightarrow Spec(\OO)$ together with an isomorphism of the generic
fibre $X_{\eta}$ with $X$. Let $Y$ be the special fibre $\XX_v=\XX
\times Spec(k(v))$. We will always also make the assumption that the
model is strictly semi-stable, which means that it is a regular
model and the fibre $Y$ is a divisor with normal crossings whose
irreducible components are smooth, have multiplicity one and
intersect transversally.

Let $i:Y \hookrightarrow \XX$ denote the inclusion map.

\subsection{Consani's Double Complex}

In \cite{cons}, Consani defined a double complex of Chow groups of
the components of the special fibre with a monodromy operator $N$
following the work of Steenbrink \cite{stee} and
Bloch-Gillet-Soul\'{e} \cite{bgs}. Using this complex she was able
to relate the higher Chow group of the special fibre at a
semi-stable prime to the regular Chow groups of the components. This
relation is what is used in defining the group $PCH$.

Let $Y=\bigcup_{i=1}^{t} Y_i$ be the special fibre of dim $n$ with
$Y_i$ its irreducible components. For $I \subset \{1,\dots,t\}$,
define
$$Y_I= \cap_{i \in I} Y_i$$
Let $r=|I|$ denote the cardinality of $I$. Define
$$Y^{(r)}:=\begin{cases} \XX & \text{ if } r=0 \\ \coprod_{|I|=r} Y_{I}& \text{ if } 1 \leq r \leq n \\
\emptyset & \text { if } r>n \end{cases}$$
For $u$ and $t$ with $1 \leq u \leq t < r$ define the map
$$\delta(u):Y^{(t+1)} \rightarrow Y^{(t)}$$
as       follows.        Let       $I=(i_1,\dots,i_{t+1})$ with
$i_1<i_2<...<i_{t+1}$. Let $J=I-\{i_u\}$. This gives an embedding
$Y_I \rightarrow Y_J$. Putting these  together induces the map
$\delta(u)$. Let $\delta(u)_*$  and $\delta(u)^*$ denote the
corresponding maps on Chow homology  and cohomology respectively.
They further  induce the Gysin and restriction maps on the Chow
groups.

Define
$$\gamma:=\sum_{u=1}^{r+1} (-1)^{u-1} \delta(u)_*$$
and
$$ \rho:=\sum_{u=1}^{r+1} (-1)^{u-1} \delta(u)^*$$
These maps have the properties that
\begin{itemize}
\item $\gamma^2=0$ \item $\rho^2=0$ \item $\gamma \cdot \rho +
\rho \cdot \gamma =0$
\end{itemize}

\subsection{The group PCH} Let $a,q$ be two integers with $q-2a >0$.
$$PCH^{q-a-1}(Y,q-2a-1):= \begin{cases} \frac{Ker(i^*i_*:CH_{n-a}(Y^{(1)})
\rightarrow CH^{a+1}(Y^{(1)}))}{Im(\gamma:CH_{n-a}(Y^{(2)})
\rightarrow CH_{n-a}(Y^{(1)}))}\otimes \Q & \text{ if $q-2a=1$ }\\
\\\frac{Ker(\gamma:CH_{n-(q-a-1)}(Y^{(q-2a)}) \rightarrow CH_{n-(q-a-1)}(Y^{(q-2a-1)}))}
{Im(\gamma:CH_{n-(q-a-1)}(Y^{(q-2a+1)}) \rightarrow
CH_{n-(q-a-1)}(Y^{(q-2a)}))}\otimes \Q  & \text { if $q-2a>1$}
\end{cases}
$$
Here $n$ is the dimension of $Y$. Note that if $q-2a>1$ and $Y$ is
non-singular, this group is $0$, while if $Y$ is singular and
semi-stable, the Parshin-Soul\'{e} conjecture implies that this
group is $CH^{q-a-1}(Y,q-2a-1)\otimes \Q$. If $q-2a=1$ and $Y$ is
non-singular, the group is $CH^a(Y) \otimes \Q$. Our interest is in
the remaining case, namely when $q-2a=1$ and $Y$ is singular.

The `Real' Deligne cohomology has the property that its dimension is
the order of the pole  of the Archimedean factor of the $L$-function
at a certain point on the left of the critical point. The group
$PCH^1(Y)$ has a similar property. Let $F^*$ be the geometric
Frobenius and $N(v)$ the number of elements of $k(v)$. The local
$L$-factor of the $(q-1)^{st}$-cohomology group is then
$$L_v(H^{q-1}(X),s)=(det(\textsl{I}-F^*N(v)^{-s}|H^{q-1}(\bar{X},\QL)^{I}))^{-1}$$
\begin{thm}[Consani] Let $v$ be a place of semistable reduction.
Assuming the weight-monodromy conjecture, the Tate conjecture for
the components   and the injectivity of the cycle class map on the
components $Y_I$, the Parshin-Soul\'{e} conjecture and that $F^*$
acts semisimply on $H^*(\bar{X},\QL)^{I}$. we have
$$\dim_{\Q} PCH^{q-a-1}(Y,q-2a-1)=-\ord_{s=a} L_v(H^{q-1}(X),s):=d_v$$
\end{thm}

\begin{proof} \cite{cons}, Cor 3.6.
\end{proof}
From this point of view the group $PCH^{q-a-1}(Y,q-2a-1)$ can be
viewed as a non-Archimedean analogue of the `Real' Deligne
cohomology. Since the $L$-factor at a prime of good reduction does
not have a pole at $s=a$ when $q-2a>1$, the Parshin-Soul\'{e}
conjecture can be interpreted as the statement that this
non-Archimedean Deligne cohomology has the correct dimension, namely
$0$, even at a prime of good reduction.

As is clear from the definition, the group $PCH$ depends on the
choice of the semi-stable model of $X$. However, Consani's theorem
says that the dimension does not. So to a large extent one can work
with any semi-stable model. Perhaps the correct definition is one
obtained by taking a limit of semi-stable models as in the work of
Bloch,Gillete and Soul\'{e} \cite{bgs} on non-Archimedean Arakelov
theory.

From this point on we specialize to the case when $X$ is a surface
and further $n=2$,$q=3$ and $a=1$. We will be interested in group
$CH^2(X,1)$ and the map to $PCH^1(Y):=PCH^1(Y,0)$. This is related
to the order of the pole of the $L$-function of $H^2(X)$ at $s=1$.
Soon we will further specialize to the case when $X=E_1 \times E_2$.

\subsection{Elements of the higher chow group}

The group  $CH^2(X,1)$ has the following presentation \cite{rama}.
It is generated by formal sums of the type
$$\sum_i (C_i,f_i)$$
where $C_i$ are curves on $X$ and $f_i$ are $\bar{K}$-valued
functions on the $C_i$ satisfying the cocycle condition
$$\sum_i \div{f_i}=0.$$
Relations in this group are give by the tame symbol of pairs of
functions on $X$.

There are some decomposable elements of this group coming from the
product structure
$$CH^1(X) \otimes CH^1(X,1) \longrightarrow CH^2(X,1)$$
A theorem of Bloch \cite{bloc} says that $CH^1(X,1)$ is simply
$K^{*}$ where $K$ is the field of definition of $X$ so such an
element looks like a sum of elements of the type $(C,a)$ where $C$
is a curve on $X$ and $a$ is in $K^*$. More generally, an element is
said to be decomposable if it can be written as a sum of products as
above over possibly an extension of the base field. Elements which
are not decomposable are said to be indecomposable.

The group $CH^2(X,1)\otimes \Q$ is the same as the $\K$-cohomology
group $H^1_{Zar}(X,\K_{2})\otimes \Q$ and the motivic cohomology
group  $H^3_{\M}(X,\Q(2))$.

\subsection{The boundary map}

The usual Beilinson regulator maps the higher chow group to the Real
Deligne cohomology. In the non-Archimedean context, it appears that
the boundary map
$$\partial:CH^2(X,1) \longrightarrow PCH^1(Y)$$
plays a similar role. It is defined as follows
$$\partial( \sum_i (C_i,f_i))=\sum_i \div_{\bar{C_i}}(f_i)$$
where $\bar{C_i}$ denotes the closure of $C_i$ in the semi-stable
model $\XX$ of $X$. By the cocycle condition, the `horizontal
divisors' namely, the closure of $\sum_i \div_{C_i}(f_i)$ cancel out
and the result is supported on the special fibre. Further, since the
boundary $\partial$ of an element is the sum of divisors of
functions, it lies in $Ker(i^*i_{*})$.

For a decomposable element of the form $(C,a)$ the regulator map is
particularly simple to compute,
$$\partial((C,a))=\ord_v(a) C_v$$

\section{Products of Elliptic Curves}

From now on we specialize to the case when $X=E_1 \times E_2$ where
$E_1$ and $E_2$ are elliptic curves over a local field $K$ of
residue characteristic $v$ with semi-stable reduction at $v$. Let
$\E_1$ and $\E_2$ denote the N\'{e}ron minimal models of $E_1$ and
$E_2$ over $S=Spec(\OO)$ respectively. The special fibre at $v$ of
the $E_i$ are N\'{e}ron polygons -
$$\E_{i,v}=\cup_{j=0}^{k_i-1}
\E_{i,v}^j$$
where  $k_i$ denotes the number of components of the special fibre
of $\E_i$. Each $\E_{i,v}^j\simeq \CP^1$. Let $\E_{i,v}^0$ denote
the identity component - that is the component which intersects
the $0$-section.

\subsection{Semi-stable models of elliptic curves}
\label{semistable}

In this section we describe the semi-stable model of the product of
elliptic curves. The product of  semi-stable models of $E_1$ and
$E_2$ is unfortunately not semi-stable - one has to blow up certain
points lying on the intersection of the products of the components .
Locally, one has the following description \cite{cons2}:

\begin{lem} Let $z_1z_2=w_1w_2$ be a local description of $\E_1
\times_{S} \E_2$ around the point $(P,Q)$ where $P$ and $Q$ are
double points lying on the intersection of two components of the
special fibre of $\E_{i,v}$, say $P \in \E_{1,v}^{0} \cap
\E_{1,v}^{1}$ and $Q \in \E_{2,v}^{0} \cap \E_{2,v}^{1}$. After a
blow up of $\E_1 \times \E_2$ with center at the origin
$(z_1,z_2,w_1,w_2)$ the resulting degeneration is normal crossings.
The special fibre $Y$ is the union of five irreducible components
$Y=\cup_{i=1}^{5} Y_i$. We label them  as follows:
$Y_1=\widetilde{(\E_{1,v}^0 \times
\E_{2,v}^0)},Y_2=\widetilde{(\E_{1,v}^0 \times \E_{2,v}^1)},
Y_3=\widetilde{(\E_{1,v}^1 \times \E_{2,v}^0)},
Y_4=\widetilde{(\E_{1,v}^1 \times \E_{2,v}^1)}$, where $\sim$
denotes the strict transform of the $(\E_{1,v}^i \times
\E_{2,v}^j)$, and $Y_5$ is the exceptional divisor. $Y_5$ is
isomorphic to $\CP^1 \times \CP^1$.
\end{lem}

\begin{proof} \cite{cons2}- Lemma 4.1.
\end{proof}

To get the global situation, we have to repeat this construction at
every double point. Let $\Z$ denote the semi-stable model and
$\psi:\Z  \rightarrow \E_1 \times \E_2$ the blowing up map.

The labeling of the components $Y_i$  is important and is with
respect to the point being blown up - each component $\E_{1,v}^{a}
\times \E_{2,v}^{b}$ is $Y_1$ with respect to the south-west corner
being blown up, $Y_2$ with respect to the north-west, $Y_3$ with
respect to the south-east and $Y_4$ with respect to the north-east.

We use this description to compute the group $PCH^1(Y)$. The group
in question is
$$PCH^1(Y)=\frac{Ker(i^*i_*:CH_{1}(Y^{(1)})
\rightarrow CH^{2}(Y^{(1)}))}{Im(\gamma:CH_{1}(Y^{(2)}) \rightarrow
CH_{1}(Y^{(1)}))} \otimes \Q$$
$Y^{(1)}$ consists of the disjoint union of the components $Y_i$ and
$Y^{(2)}$ consists of their pairwise intersections, $Y_i \cap Y_j$,
denoted by $Y_{ij}$. Similarly, the intersections of three
components $Y_i \cap Y_j \cap Y_k$ will be denoted by $Y_{ijk}$.

From the description above we see that $Y_{ij}$ is one of
`horizontal curves'  $Y_{12}$ and $Y_{34}$ which are $\E_{1,v}^{0}
\times Q$ and $\E_{1,v}^{1} \times Q$ respectively, `vertical
curves' $Y_{13}$ and $Y_{24}$ which are $P \times \E_{2,v}^{0}$ and
$P \times \E_{2,v}^{1}$  or `exceptional curves' $Y_{i5},
i=\{1,\dots 4\}$. All the curves $Y_{ij}$ are isomorphic to $\CP^1$.

\begin{prop} When $X=E_1 \times E_2$ the group $PCH^1(Y)$ is three
dimensional. It is  generated by $\mE_1=\psi^*(\E_{1,v} \times Q)$,
$\mE_2=\psi^*(P \times \E_{2,v})$, which are the images of $\E_1$
and $\E_2$ in the special fibre of the semi-stable model, and the
sum over all the exceptional divisors  of the cycle $Y_{15} + Y_{45}
- Y_{25} -Y_{35}$, $\mF= \sum (Y_{15} + Y_{45} - Y_{25} -Y_{35})$.
\end{prop}
\begin{proof}

We first show that our cycles lie in the kernel of $i^*i_*$.

\begin{lem} The elements $\mF$,$\mE_1$ and $\mE_2$ lie in the $Ker(i^{*}i_{*})$.
\end{lem}

\begin{proof} The group which is the target of $i^*i_*$ is
$$CH^2(Y^{(1)})=\bigoplus_j CH^2(Y_j)$$
where the sum is over all the components $Y_j$ which arise from all
the blow ups. A cycle lies in the kernel of $i^{*}i_{*}$ if the
restriction of its image under $i_*$ to the $CH^2(Y_j)$ is $0$ for
all $j$.

Around a the blow up of a point, in the group $\bigoplus_{j=1}^{5}
CH^2(Y_i)$ the image of $i^*i_*(Y_{15} + Y_{45} - Y_{25} -Y_{35})$
is
$$(-Y_{125}-Y_{135},Y_{125}+Y_{245},Y_{135}+Y_{345},-Y_{245}-Y_{345},0)$$
As remarked above, each $Y_i$ is $Y_1$, $Y_2$,$Y_3$ and $Y_4$
depending on which corner is being blown up, so the image of the
cycle $\mF$ under $i^*i_*$ in each $Y_i$, $i\neq 5$ is
$$(-Y_{125}-Y_{135}+Y'_{125}+Y'_{245}+Y''_{135}+Y''_{345}-Y'''_{245}-Y'''_{345}).$$
This is rationally equivalent to $0$ as any $0$-cycle of degree $0$
on $\CP^1 \times \CP^1$ is rationally equivalent to $0$.

The cycles $\mE_1$ and $\mE_2$ lie in $Ker(i_*i^*)$ as they are the
restrictions of $div(\pi)$, where $\pi$ is the uniformizer at $v$,
to the cycles $\E_1$ and $\E_2$ respectively. Consani shows that
$\psi^*(\E_{1,v})$ in $\bigoplus_{j=1}^5 CH^1(Y_j)$ is
$Y_{12}+Y_{34}+(Y_{15}-Y_{45}+Y_{25}-Y_{35})$. As we will see below,
the cycles supported in $Y_5$ are rationally equivalent to $0$ so
$\mE_1$ is simply the closure of $\E_{1,v} \times Q$ in $PCH^1(Y)$.
The same holds for $\mE_2$.

\end{proof}

The cycles $Y_{15}$ and $Y_{45}$ do not intersect as the components
$Y_1$ and $Y_4$ do not intersect. Further, they are reduced. Hence
these two cycles are parallel and hence rationally equivalent in
$CH^1(Y_5)$. Similarly $Y_{25}$ and $Y_{35}$  are rationally
equivalent in $CH^1(Y_5)$. The cycles $Y_{15}$ and $Y_{25}$
intersect precisely at one point, $Y_{125}$, with multiplicity one,
so they give rulings of $Y_5 \simeq \CP^1 \times \CP^1$. Further,
the cycle  $Y_{15}$ lying in the intersection of a $\CP^1 \times
\CP^1$ with the exceptional fibre of the blow up of the south west
corner, is equivalent to the cycle $Y_{45}'$ which is  the
intersection of the same $\CP^1 \times \CP^1$ with the exceptional
fibre of the blow up of the north east corner. Similarly for
$Y_{25}$ and $Y_{35}$. Hence the cycles $Y_{15}-Y_{25}$ lying in all
the exceptional divisors are all equivalent in the group $PCH^1(Y)$.
So $\mF=(k_1 k_2)(Y_{15}+Y_{45}-Y_{25}-Y_{35})$, where $k_1$ and
$k_2$ are the number of components of $\E_{1,v}$ and $\E_{2,v}$
respectively.

Similarly, the cycles $\E_{1,v}^{j} \times P$ and $\E_{1,v}^j \times
P'$ in $\E_{1,v}^{j} \times \E_{2,v}^{l}$ are equivalent for any two
points  $P$ and $P'$. Further, the cycles $Y_{ij}$ embedded in $Y_i$
and $Y_j$ are equivalent in the group $PCH^1(Y)$ as their difference
lies in the image of the Gysin map $\gamma$. As a result, for a
fixed $j$ any two  cycles of the form $\E_{1,v}^{j} \times Q $ for a
$Q$ in $\E_{2,v}$ are equivalent, and similarly, for a fixed $l$ ,
any two cycles of the form $P \times \E_{2,v}^{l}$ for $P$ in
$\E_{1,v}$ are equivalent under the image of the Gysin map.

So the result is that we have three cycles in the group $PCH^1(Y)$,
namely  $\mE_1$, $\mE_2$ and $\mF$ which are clearly linearly
independent.

As one can compute the order of the pole of the local L-factor, the
theorem of Consani \cite{cons} asserts that this group is three
dimensional, so these cycles generate the group.

\end{proof}

To prove surjectivity, therefore,  we have to find three elements of
the Chow group $CH^2(E_1 \times E_2,1) \otimes \Q$ which bound these
three generators.

\begin{rem} When $E_1 \times E_2$ are elliptic curves over
$\Q$ the dimension of the Real Deligne cohomology, which is the
target of the Beilinson regulator map, is also three dimensional. In
that case it is easy to find cycles which bound two of the three
generators. The third requires more work - one has to use the
modular parametrization \cite{beil}.
\end{rem}

\begin{rem} At a prime of good reduction for both $E_1$ and $E_2$,
in the group $PCH^1(Y)=CH^1(Y) \otimes \Q$, one has the cycles
$\mE_1=\E_{1,v}$ and $\mE_2=\E_{2,v}$ as well as cycles for any
isogeny between $\E_{1,v}$ and $\E_{2,v}$. So there are either
$2$,$4$ or $6$ cycles depending on whether $\E_{1,v}$ and $\E_{2,v}$
are not isogenous, isogenous, but not supersingular, or have
supersingular reduction.
\end{rem}

\subsection{Genus two curves on products of elliptic curves}

Spei{\ss} \cite{spie} constructed an element of the higher chow
group using a genus two curve on the generic fibre. We show that his
construction can be used in our case of semistable reduction as
well. We have to use some work of Frey and Kani \cite{frka} on the
existence of irreducible genus two curves whose Jacobian is
isogenous to a product of elliptic curves.

\begin{thm}[Frey and Kani]
Let $K$ be a local field with residual characteristic $v$ and $E_1$
and $E_2$ two elliptic curves over $K$. Let $n$ be an odd integer,
$E_1[n]$ and $E_2[n]$ the $n$-torsion subgroups and $\phi:E_1[n]
\rightarrow E_2[n]$ a $K$-rational homomorphism which is an
anti-isometry with respect to the Weil pairings - that is
$e_n(\phi(x),\phi(y))=e_n(x,y)^{-1}$. Let $J=E_1 \times
E_2/graph(\phi)$ and $p:E_1 \times E_2 \rightarrow J$ the
projection. Then there exists a unique curve $C \subset J$ defined
over $K$ such that the following holds.

\begin{itemize}
\item $C$ is a stable curve of genus two in the
sense of Deligne and Mumford.

\item $-id^*(C)=C$.

\item Let $\lambda_{C}$ denote the map from $J \rightarrow
\check{J}$ induced by the line bundle corresponding to $C$. Then the
composite maps,
$$\pi_i:C \stackrel{j}{\longrightarrow} J \stackrel{\lambda_{C}}{\longrightarrow}
\check{J} \stackrel{\check{p}}{\longrightarrow}E_1 \times E_2
\longrightarrow E_i \text{ i=1,2}$$
are finite morphisms of degree $n$.
\end{itemize}
\end{thm}
\begin{proof}\cite{frka}, Proposition [1.3].
\end{proof}

We now apply this criterion in a special case, choosing $n$
judiciously so as to ensure that we bound the right cycle. This is
a variation of the method used in \cite{spie}, Lemma [3.3].

Let $a$ be an integer and $n=a^2+1$. Choose $a$ such that
$(a,v)=(n,v)=1$ and $n$ is odd. Extend $K$ to a field where all the
$n$-torsion of $E_1$ and $E_2$ are defined. From the theory of
N\'{e}ron models, we have that $n$ then divides the number of
components $k_i$ of the special fibres  $\E_{i,v}$ of the N\'{e}ron
models $\E_i$ of $E_i$.

As a group the special fibre $\E_{i,v}$ is isomorphic to $\G_{m}
\times \ZZ/k_i\ZZ$. We will denote an element of $\E_{i,v}$ by
$(x;m)$ with $x \in \G_m$ and $m \in \ZZ/k_i\ZZ$. Let $h_a$ denote
the isogeny
$$h_a:\E_{1,v} \longrightarrow \E_{2,v}$$
$$h_a((x;m))=(x^a;a \cdot m)$$
where multiplication by $a$ is to be understood as the action of the
class of $a$ in $\ZZ/(k_1,k_2)\ZZ$ which is identified with
$Hom(\ZZ/k_1\ZZ,\ZZ/k_2\ZZ)$. Explicitly, this is as follows. If $m$
is in $\ZZ/k_1\ZZ$ and $a$ is in $\ZZ/(k_1,k_2)\ZZ$ then
$$a \cdot m = am k_2/(k_1,k_2)  \text{ mod } k_2$$
Since $n|(k_1,k_2)$ the group $\ZZ/(k_1,k_2)\ZZ$  is non-trivial.

Let $h_a[n]:\E_{1,v}[n] \longrightarrow \E_{2,v}[n]$ denote the
restriction of $h_a$ to the $n$-torsion points. Since $(n,v)=1$ the
groups $E_i[n]$ and $\E_{i,v}[n]$ are isomorphic \cite{seta}. So the
map $h_a[n]$ lifts to a map $\phi_a:E_1[n] \rightarrow E_2[n]$.

\begin{lem} The map $\phi_a$  is an anti-isometry with respect to the Weil pairing $e_n$
\end{lem}
\begin{proof} If $X$ and $Y$ are two points in $E_1[n]$ mapping to $(x;m_x)$ and
$(y;m_y)$ in $\E_{1,v}[n]$ respectively we have
$$e_n(\phi_a(X),\phi_a(Y))=e_n(h_a((x;m_x)),h_a((y;m_y)))=e_n((x;m_x),\check{h}_a
\circ h_a((y;m_y)))$$
as the dual isogeny $\check{h}_a$ is the adjoint of $h_a$ with
respect to the Weil pairing. So this is equal to
$$e_n((x;m_x),(y^{a^2};a\cdot a \cdot m_y))=e_n((x;m_x),(y^{n-1};(n-1) \cdot m_y))=
e_n((x;m_x),(y^{-1};-1 \cdot m_y))=e_n(x,y)^{-1}$$
as  $\check{h}_a \circ h_a$ is multiplication by $a^2=n-1$ and $Y$
is in the $n$ torsion, so $(n-1) \cdot m_y=-1 \cdot m_y$ and
$y^{n-1}=y^{-1}$.

\end{proof}

From the theorem of Frey and Kani with $\phi=\phi_a$ we get a
corresponding stable genus $2$ curve $C$ and finite morphisms
$\pi_i:C \rightarrow E_i$ of degree $n$. $C$ is a principal
polarization on $J=E_1 \times E_2/(graph (\phi))$. It satisfies the
additional property that $p^*(C) \sim n\Theta$, where $\Theta=E_1
\times 0 \cup 0 \times E_2$. Further, it is the unique curve
satisfying that as well as $-id^{*}(C)=C$.

We would like to understand the special fibre of the closure of this
curve in a semistable model of $E_1 \times E_2$. We first describe
what happens in the product of the two N\'{e}ron models $\E_1 \times
\E_2$ and then describe its image in the semi-stable model of $E_1
\times E_2$ constructed in section \ref{semistable}. Let $\CC$
denote the closure of $C$ in the N\'{e}ron model of $J$, $\J$.

\begin{prop} The special fibre $\CC_v$ of  $\CC$ is reducible and
is isomorphic to $\E_{1,v} {\sqcup}_{1} \E_{2,v}$, namely a union of
two curves isomorphic to $\E_{1,v}$ and $\E_{2,v}$ which meet
transversally at $((1;0),(1;0))$, the identity, in the product of
$\E_{1,v}^{0}$ and $\E_{2,v}^{0}$, the product of the identity
components, and nowhere else. The finite maps $\tilde{\pi}_i:\CC_{v}
\longrightarrow \E_{i,v}$ are given by $\tilde{\pi_1}=id
{\sqcup}{_1} \check{h}_{-a}$ and $\tilde{\pi_2}=h_a {\sqcup}_{1}
id$.

\label{specialfibre}
\end{prop}

\begin{proof} We follow the argument in Spie{\ss}\cite{spie} mutatis
mutandis. Let $\Theta_v=\E_{1,v} \times (1;0) \cup (1;0) \times
\E_{2,v}$. This is a stable genus $2$ curve on $\E_{1,v} \times
\E_{2,v}$. The idea is to show that $\Theta_v$ satisfies the
properties of Frey-Kani \cite{frka}, Proposition 1.1, which states
that there is a unique stable curve of genus two $\CC_v'$ which
satisfies the conditions that
\begin{itemize}
\item $-id^*(\CC'_v)=\CC'_v$
\item $p^*(\CC'_v)=n\Theta_v$, where $p$ is as described below.
\end{itemize}
Since the special fibre $\CC_v$ satisfies these, and we will show
that $\Theta_v$ satisfies these, we will have that $\Theta_v\simeq
\CC_v\simeq \CC'_v$.

Clearly, $-id^*(\Theta_v)=\Theta_v$. To show the second property, we
proceed as follows. Let
$$p=\begin{pmatrix} id & \check{h}_a\\ h_{-a} & id
\end{pmatrix}:\E_{1,v} \times \E_{2,v} \rightarrow \E_{1,v} \times \E_{2,v}$$
Let $X=(x;m_x)$ be an element of $\E_{1,v}[n]$. Then
$$p(X,h_a[n](X))=(X.\check{h}_a \circ
h_a(X),h_{-a}(X).h_a(X))$$
$$=((x^{1+a^2};(1+a^2) \cdot m_x),(x^{-a+a},(-a+a) \cdot m_x))=((1;0),(1;0))$$
as $a^2+1=n$ and $X$ is in the $n$-torsion. So $graph(h_a[n])
\subset Ker(p)$. Similarly, we can see that the kernel of
$$\check{p}=\begin{pmatrix} id & \check{h}_{-a}\\ h_a & id
\end{pmatrix}:\E_{1,v} \times \E_{2,v} \rightarrow \E_{1,v} \times \E_{2,v}$$
contains $graph(\check{h}_a[n])$. Since $\check{p} \circ p=[n]$ we
have
$$|Ker(\check{p} \circ
p)|=|Ker(p)|^2=|Ker([n])|=n^4=|graph(h_a[n])|^2$$
so $Ker(p)=graph(h_a[n])$ and we can identify the image of $p$ with
$\E_{1,v} \times \E_{2,v} /(graph (h_a[n])$.
\begin{lem}  $p^*(\Theta_v)\sim n\Theta_v$.
\end{lem}

\begin{proof}

Let $\Gamma_{h_a}$, respectively  $\Gamma_{\check{h}_{-a}}^t$ denote
the graphs of the maps $(id,h_a)$, respectively
$(\check{h}_{-a},id)$ from $\E_{1,v}$, respectively $\E_{2,v}$, to
$\E_{1,v} \times \E_{2,v}$. Since the diagrams

\begin{center}
\begin{tabular}{cp{1in}c}
$
\begin{CD}
\E_{1,v} @>(id,h_a)>> \E_{1,v} \times \E_{2,v}\\
@V[n]VV @VVpV \\
\E_{1,v} @>(id,0)>>  \E_{1,v} \times \E_{2,v}
\end{CD}
$
&&
$
\begin{CD}
\E_{2,v} ,\, @>(\check{h}_{-a},id)>> \E_{1,v} \times \E_{2,v} \\
@V[n]VV @VVpV \\
\E_{2,v} @>(id,0)>> \E_{1,v} \times \E_{2,v}
\end{CD}
$
\end{tabular}
\end{center}
are cartesian, we have
$$p^*(\Theta_v)=\Gamma_{h_a} \cup \Gamma^t_{\check{h}_{-a}}$$
Since the closure of $\E_{1,v} \times \E_{2,v}$ is a union of $\CP^1
\times \CP^1$'s the divisor $p^*(\Theta )$ can be written as a sum
$p_1^*(D_1) + p_2^*(D_2)$ where $p_1$ and $p_2$ are the projection
maps to $\E_{1,v}$ and $\E_{2,v}$ respectively and $D_1$ and $D_2$
are divisors on $\CP^1$.

We have
$$[\Gamma_{h_a}].[\E_{1,v} \times (1;0)]=[Ker(h_a) \times (1;0)]$$
$$[\Gamma_{h_a}].[(1;0) \times \E_{2,v}]=[(1;0) \times
(1;0)]=[\Gamma_{\check{h}_{-a}}^t].[\E_{1,v} \times (1;0)]$$
$$[\Gamma_{h_{-a}^t}].[(1;0) \times Ker(\check{h}_{-a})]$$
From this we get,
$$D_1 \sim (p_1)_*(p_1^*(D_1)+p_2^*(D_2))=[Ker(h_a)]+[(1;0)]$$
and similarly $D_2 \sim [(1;0)]+ [Ker(\check{h}_{-a})]$. Since any
two points on $\CP^1$ are equivalent, we have $D_1$ and $D_2$ are
equivalent to $(a^2+1)[(1;0)]=n[(1;0)]$. So we have
$$p^*(\Theta_v)=p_1^*(n[(1;0)])+p_2^*(n[(1;0)])=n(\E_{2,v}+\E_{1,v})=n \Theta_v$$

\end{proof}

The map $\lambda_C$ and $\lambda_{\Theta}$, which are isomorphisms
from $J \rightarrow \check{J}$ and $E_1 \times E_2 \rightarrow
\check{(E_1 \times E_2)}$ induced by the principal polarizations $C$
and $\Theta$ extend to isomorphisms of the N\'{e}ron models and in
particular, induce isomorphisms of the special fibres. The map $p$
induces a homomorphism $\check{p}:\check{J} \rightarrow \check{(E_1
\times E_2)}$ which is the same as $p^*$ on the divisors of degree
$0$.

Let $p'$ be the homomorphism
$$p'=\lambda_{\Theta}^{-1} \circ
\check{p} \circ \lambda_{C}$$
extended to induce a homomorphism of the special fibres. From the
definition, it is easy to see $Ker(p')$ is contained in the
$n$-torsion.

If ${\CC'}_v$ is a stable genus $2$ curve satisfying
$p^*({\CC'}_v)=n\Theta_v$ then one has ${\CC'}_v=T_x(\CC_v)$ for
some $x$ in $Ker(p')$. As $n$ is odd, if ${\CC'}_v$ further
satisfies the condition that $-id^*({\CC'}_v)={\CC'}_v$ then
${\CC'}_v \simeq \CC_v$, otherwise it would imply that there is an
element of $2$-torsion in $Ker(p')$. Since $\Theta_v$ satisfies this
additional condition, $\Theta_v \simeq \CC_v$. The rest of
Proposition \ref{specialfibre} follows by chasing the definitions of
the various maps.

\end{proof}

\subsection{A new element}

Using the genus 2 curve constructed above we can get a new element
of $CH^2(E_1 \times E_2,1)$ - making a clever  choice of a pair of
Weierstra{\ss } points on the genus two curve.

The construction is as follows. Let $\CC'$ be a minimal regular
model of $C$. From Parshin \cite{pars} we have a description of the
special fibre as well as a description of the closure of the
Weierstra{\ss} points on the special fibre.

In our case the special fibre has the following description (VI, in
Parshin's notation) - there are two genus 0 curves, $B_1$ and $B_2$
with self intersection $-3$. To each of these is attached a chain of
genus $0$ curves $X_i, i=\{1,\dots,r\}$ and $Z_k,k=\{1\dots
t\}$,with $t$ and $r$ odd, respectively with self intersection $-2$,
such that each curve intersects the neighboring two curves at a
single point.

In other words, these are the N\'{e}ron special fibres of semistable
elliptic curves. The two semi-stable fibres of elliptic curves are
joined by a chain of genus $0$ curves $L_j,j=\{1,\dots,s\}$ with
self intersection $-2$ which meet at the identity components. So in
particular, $r=k_1-1$, $t=k_2-1$ and $B_1$ and $B_2$ correspond to
the identity components of $\E_{1,v}$ and $\E_{2,v}$ respectively.

The closure of the Weiersta{\ss} points is as follows - one point
lies on each $B_1$ and $B_2$ and two points each intersect the
components $X_{\frac{s+1}{2}}$ and $Z_{\frac{t+1}{2}}$.

In particular, we have a function $f_{P,Q}$ on $C$ such that the
closure of $P$ lies on $B_1$ and the closure of $Q$ lies on $B_2$
and whose divisor on $C$ is
$$\div(f_{P,Q})=2(P)-2(Q).$$
The divisor of $f_{P,Q}$ on $\CC'$ can be expressed in terms of the
components above
$$\div_{\CC'}(f_{P,Q})={\mathcal H} + a_1B_1 + \sum_{i=1}^r b_i
X_i + \sum_{j=1}^{s}c_j L_j + \sum_{k=1}^{t} d_k Z_k + a_2 B_2$$
where ${\mathcal H}$ is the closure of the divisor $2(P)-2(Q)$ -
that is, the horizonal component. Multiplying $f_{P,Q}$ by a power
of the uniformizer $\pi$ one can assume that $a_2=0$ as
$\div_{\CC'}(\pi)=\CC'_v$.
\begin{lem} If $f_{P,Q}$ is as above with $a_2=0$ then $a_1\neq 0$.
\label{xi}
\end{lem}
\begin{proof}
Since $\CC'$ is a minimal regular model we can use the intersection
theory of arithmetic surfaces described, for example, in Lang
\cite{lang}, Chapter 3. In particular, we have that the intersection
number
$$(\div_{\CC'}(f_{P,Q}).D)=0$$
for {\em any} divisor $D$ contained in the special fibre. Applying
this to different choices of $D$, namely $D=B_i,X_i,L_j,Z_k$ and
using what we know of their intersections and self-intersections
gives us the following set of equations -
$$-3a_1+b_1+b_r+c_1+2=0$$
$$a_1-2b_1+b_2=0$$
$$b_{i-1}-2b_{i}+b_{i+1}=0 \;\;\;\; \{i=2\dots r-1\}$$
$$b_{r-1}-2b_r+a_1=0$$
$$a_1-2c_1+c_2=0$$
$$c_{j-1}-2c_j+c{j+1}=0 \;\;\;\; \{j=2 \dots s-1\}$$
$$c_{s-1}-2c_s=0$$
$$-2d_1+d_2=0$$
$$d_{k-1}-2d_k+d_{k+1}=0 \;\;\;\; \{k=2\dots t-1\}$$
$$d_{t-1}-2d_t=0$$
$$c_s+d_1+d_t-2=0$$
Solving these equations shows $d_k=0,k=\{1,\dots,t\}$,
$c_j=2(s+1-j)$, so in particular $c_s=2$,$c_1=2s$ and finally
$a_1=b_i, i=\{1 \dots r\}=2(s+1)$. In particular, since $s\geq 0$ we
have $a_1 \neq 0$.

\end{proof}

Recall that we have maps $\pi_i:C \rightarrow E_i$, which induce a
map $\rho:C \rightarrow E_1 \times E_2$. $\rho$ is generically a
closed immersion. Let $P_i=\pi_i(P)$ and $Q_i=\pi_i(Q)$. There are
functions $f_1$ on $E_1 \times P_2$ and $f_2$ on $Q_1 \times E_2$
with
$$\div(f_1)=2(P_1,P_2)-2(Q_1,P_2) \text{ and }
\div(f_2)=2(Q_1,P_2)-2(Q_1,Q_2)$$
On the closure, by the description of the maps $\pi_i$ on the
special fibre, both $P$ and $Q$ map to the identity components of
the special fibres $\E_i^0$ of $E_i$. Hence the divisors of $f_i$
in the semistable model of $E_1 \times E_2$ do not contain any
components of the special fibre.

Define $\Xi=\Xi_{P,Q}$ be the cycle
$$\Xi=(C,f_{P,Q})+(E_1 \times P_2,f_1^{-1})+(Q_1 \times
E_2,f_2^{-1})$$
From the definition of $P_i$ and $Q_i$ we have
$$\div_{C}(f_{P,Q}))-\div(f_1)-\div(f_2)$$
$$=2(P_1,P_2)-2(Q_1,Q_2)-2(P_1,P_2)+2(Q_1,P_2)-2(Q_1,P_2)+2(Q_1,Q_2)=0$$
hence $\Xi$ is an element of $CH^2(E_1 \times E_2,1)$.

Since $\div(f_i)$ have no components in the special fibre,
$$\div_{\CC'}(f_{P,Q})+ \div_{\E_1 \times  P_2}(f_1^{-1}) + \div_{Q_1 \times \E_2}(f_2^{-1})  =
(2s+2)\left(B_1+\sum_{i=1}^{r} X_i\right) + \sum_{j=1}^{s} c_j L_j$$
where $c_j$ are as above. In the next section, we relate this to the
cycles in $PCH^1(Y)$.

\subsection{Surjectivity of the boundary map}

In this section we compute the image of the elements under the
boundary map
$$\partial:CH^2(E_1 \times E_2,1) \rightarrow
PCH^1(Y)$$
We have two decomposable elements, $(E_1 \times 0,\pi) $ and $(0
\times E_2,\pi)$, whose boundary is
$$\partial((E_1 \times 0,\pi))=\psi^*(\E_{1,v} \times (1;0))=\mE_1$$
$$\partial((0 \times E_2,\pi))=\psi^*((1;0) \times \E_{2,v})=\mE_2$$
We also have the third element $\Xi_{P,Q}$. To compute its boundary
observe that under the map $\CC' \rightarrow \CC$ the linking
components $L_j$  of the special fibre collapse to a point. Further,
the $X_i$ and $B_1$ map on to the graph $\Gamma_{h_a}$ and
similarly, the $Z_k$ and $B_2$ map to $\Gamma_{\check{h}_{-a}}^{t}$.
So, for the choice of $f_{P,Q}$ above, the boundary of $\Xi$ in
$\E_1 \times \E_2$ is
$$\partial(\Xi)=(2s+2) \Gamma_{h_{a}}$$

The points at which the curves $X_i$ meet each other or $B_1$ are
being blown up. The graph $\Gamma_{h_a}$ lies in either the
component $Y_1$ or $Y_4$ with respect to the point being blown up.
As a consequence the total transform is
$$\psi^*(\Gamma_{h_a})=\bar{\Gamma}_{h_a}+(k_1,k_2)(Y_{15}+Y_{45})$$
as $\Gamma_{h_a}$ has $(k_1,k_2)$ components and  where
$\bar{\Gamma}_{h_a}$ denotes the closure of $\Gamma_{h_a}$ in the
blow up. Similarly, the boundary of the decomposable element
$(C,\pi)$ is
$$\bar{\Gamma}_{\check{h}_a}^{t}+\bar{\Gamma}_{h_a}+(k_1,k_2)(Y_{15}+Y_{45}+Y_{35}+Y_{25}).$$
as $\CC_v=\Gamma_{\check{h}_a}^{t}+\Gamma_{h_a}$.

Therefore, the element $\Xi_{P,Q}-(s+1)(C,\pi)$ has boundary
$$\partial(\Xi_{P,Q}-(s+1)(\CC,\pi))=(s+1)\left( \bar{\Gamma}_{\check{h}_a}^{t}-
\bar{\Gamma}_{h_a}+(k_1,k_2)(Y_{15}+Y_{45}-Y_{35}-Y_{25}) \right)$$
$$=(s+1)\left(\bar{\Gamma}_{\check{h}_a}^{t}-\bar{\Gamma}_{h_a}+\frac{(k_1,k_2)}{k_1k_2}\mF\right)$$
We have
$$\bar{\Gamma}_{h_a}=\mE_1+a\mE_2 \text{ and }
\bar{\Gamma}_{\check{h}_a}^{t}=a\mE_1+\mE_2$$
so we can further subtract $(s+1)(a-1)(\mE_1-\mE_2)$ to get an
element whose boundary is precisely
$\frac{(s+1)(k_1,k_2)}{k_1k_2}\mF$.
\begin{thm} Suppose $E_1$ and $E_2$ are two non-isogenous elliptic
curves over a local field $K$ with split multiplicative reduction at
the closed point $v$. Let $\XX$ be a semi-stable model of the
product $E_1 \times E_2$ and $\XX_v$ denote the special fibre. Then
the map
$$\partial:CH^2(E_1 \times E_2,1)\otimes \Q \longrightarrow
PCH^1(\XX_v) \otimes \Q$$
is surjective.
\end{thm}

\section{Final Remarks}

Let $\Sigma$ be the group,
$$\Sigma=Ker(CH^2(\XX) \rightarrow CH^2(X)).$$
Spie{\ss} \cite{spie}, Section 4, describes some consequence of the
assumption that it is a torsion group. Surjectivity of the map
$\partial$ implies this, in fact, it implies finiteness, so all the
consequences apply in our case.

This paper began as an attempt to prove the  $S$-integral Beilinson
conjecture when $X$ is a product of two non-isogenous modular
elliptic curves over $\Q$. This remains to be done - unfortunately,
while our and Spie{\ss}' elements  can be lifted to number fields,
they cannot be used to produce surjectivity as there may be several
primes at which their boundary is non-trivial, so at best one can
get relations between codimension $2$ cycles of the type described
in the previous paragraph.

In the case of isogenous elliptic curves, Mildenhall's elements have
a boundary at precisely one prime. It is hoped that a careful
analysis of the special fibre at a semi-stable prime combined with
Mildenhall's construction may work to prove surjectivity for the
product of non-isogenous elliptic curves over $\Q$ which have bad
semistable reduction at that prime. Curiously, though it is not
clear how one can construct enough elements to prove surjectivity at
primes $p$  of good reduction where the curves may become isogenous
mod $p$.

\bibliographystyle{alpha}
\bibliography{references}

\begin{tabular}{lp{.5in}l}
{\bf Permanent Address:} && {\bf Current Address ( valid until 31 March 2008) :}\\
Ramesh Sreekantan && Ramesh Sreekantan\\

School of Mathematics && TIFR Centre for Applicable Mathematics\\
Tata Institute of Fundamental Research && Post Bag No 3\\
Homi Bhabha Road, Colaba &&  Sharadanagar, Chikkabommasandara\\
Mumbai, 400 005 India && Bangalore 560 065 India\\
\\
{\bf email:} {\em ramesh@math.tifr.res.in} && {}
\end{tabular}

\end{document}